\documentclass{raex}
\usepackage{amsmath, amssymb, amsthm, amsfonts}

\begin{Author} 
	\FirstName{Jan A.}\LastName{Grzesik}
	\PostalAddress{Allwave Corporation, Torrance, CA 90503, U.S.A}
\end{Author}

\begin{MathReviews}
	\primary{11B68,\,65Q30,\,30E20}
	\secondary{\\42A16.}
\end{MathReviews}

\begin{KeyWords} 
\end{KeyWords}

\title{Contour Integration Underlies Fundamental Bernoulli Number Recurrence}

\begin{document}

\maketitle

\begin{abstract}
%All papers in the \rae{}  must have abstracts.
      One solution to a relatively recent American Mathematical Monthly problem [6], requesting the
evaluation of a real definite integral, could be couched in terms of a contour integral which vanishes
{\textit{a priori.}}  While the required real integral emerged on setting to zero the real part of the
contour quadrature, the obligatory, simultaneous vanishing of the imaginary part alluded to still
another pair of real integrals forming the first two entries in the infinite log-sine sequence,
known in its entirety.  It turns out that identical reasoning,
utilizing the same contour but a slightly different analytic function thereon, sufficed not only to
evaluate that sequence anew, on the basis of a vanishing real part, but also, in setting to zero its
conjugate imaginary part, to recover the fundamental Bernoulli number recurrence.  The even order Bernoulli
numbers $B_{2k}$ entering therein were revealed on the basis of their celebrated connection to Riemann's zeta
function $\zeta(2k).$  Conversely, by permitting the related Bernoulli polynomials to participate as integrand factors,
Euler's connection itself received an independent demonstration, accompanied once more by an elegant log-sine
evaluation, alternative to that already given.  
And, while the Bernoulli recurrence is intended to enjoy here the pride of place, this note ends on a gloss wherein
all the motivating real integrals are recovered yet again, and in quite elementary terms, from the Fourier
series into which the Taylor development for Log$(1-z)$ blends when its argument $z$
is restricted to the unit circle.
\end{abstract}

\section{Introduction.}

      An American Mathematical Monthly problem posed within relatively recent memory [6] sought the
evaluation
\begin{equation}
\int_{\,0}^{\,\pi/2}\left\{\rule{0mm}{4mm}\log(\,2\sin(x)\,)\right\}^{2}dx\,=\,\frac{\,\,\,\pi^{\,3}\,}{\,24\,}\,.
\end{equation}
One mode of solution depended upon integration of an analytic function around the periphery $\,\Omega\,$ of a semi-infinite
vertical strip with no singularities enclosed, the quadrature having thus a null outcome\footnote{Both
contour $\,\Omega\,,$ a vertical rectangle of unlimited height, and the notion of integrating an analytic function thereon so
as to obtain a null result, imitate a similar ploy utilized in [1, Section 5.3, example 5] on behalf of (2) and still
further attributed there to Ernst Lindel{\"{o}}f.} known in advance on the strength of Cauchy's theorem.
Evaluation (1) was automatically produced by setting to zero the real part of that integral,\footnote{Two solutions for (1) were submitted by
the undersigned, one involving contour integration in the manner suggested, and the other based upon a Fourier series.
The Bernoulli recurrence (5) was assembled as a spontaneous by-product of an ancillary, null-quadrature calculation upon
that same contour $\,\Omega\,,$ initially aimed only at evaluating the log-sine integrals (4).  This note embodies the content
of that collateral calculation, slightly rephrased so as to highlight the newly recovered Bernoulli number sum identity.}
whereas the complementary requirement that the imaginary part likewise vanish brought into play, and successfully so, both known quadratures
\begin{equation}
\int_{\,0}^{\,\pi}\log(\,\sin(x)\,)\,dx\,=\,-\,\pi\log(2)
\end{equation}
and
\begin{equation}
\int_{\,0}^{\,\pi}x\log(\,\sin(x)\,)\,dx\,=\,-\frac{\,\,\,\pi^{\,2}\,}{\,2\,}\log(2)\,.
\end{equation}
With (2) and (3) in plain view, a temptation arose to provide for them, too, an {\textit{ab initio}} verification, and,
more even than that, to evaluate the entire hierarchy of log-sine integrals\footnote{If that were the only goal then we should
assuredly stop dead in our tracks, simply because, on the one hand, {\em{Mathematica}} provides all such evaluations
on demand, with great aplomb, and this even in its symbolic mode, while, on the other, a relatively painless
derivation (Eq. (51)) can be based upon a Fourier series, one which materializes in its turn from the power series for $\,{\rm{Log}}(1-z)\,$
when argument $\,z\,$ is forced to lie upon the unit circle.  This Fourier series underlies in addition an essentially zinger
verification of (1).  All such manifold benefits of the Fourier option are sketched in an Appendix.  And, prior even to that, when
Bernoulli polynomials are admitted as integrand factors during contour integration (Section 4), the log-sine series is produced
once more in (43).  It goes without saying that contour-based (Eqs. (14) and (43)) and Fourier-based (Eq. (51))
evaluations of (4), even though they may be of secondary interest in the present context, do stand in complete agreement.} 
\begin{equation}
I_{n}\,=\,\int_{\,0}^{\,\pi}x^{n}\log(\,\sin(x)\,)\,dx
\end{equation}
as the power of $\,x\,$ roams over all non-negative integers $\,n\geq 0.$
Not only was this fresh ambition, digressive and self-indulgent though it may have been, easy to satisfy via quadrature
on the same contour as before, but it also exposed to view once more the fundamental Bernoulli number recurrence
\begin{equation}
\sum_{k\,=\,0}^{n-1}  \left(\!\!\begin{array}{c}
                                     n \\
                                     k
                                  \end{array}\!\!\right)
B_{k}\,=\,0
\end{equation}
which is valid for $\,n\geq 2\,$ and, together with the initial condition $\,B_{0}\,=\,1\,$ and the self-consistent choice
$\,B_{1}\,=\,-1/2\,,$ is adequate to populate the entire Bernoulli ladder, complete with null entries at all odd indices
beyond $\,k\,=\,1\,,$ {\textit{viz.,}} $\,B_{2l+1}\,=\,0\,$ whenever $\,l\geq1\,.$  Source material on the
Bernoulli numbers and the related Bernoulli polynomials is ubiquitous, and can be sampled, for example, in [2, 10, 12].
References [3, 7] provide a valuable overview all at once of their
mathematical properties and historical genesis in computing sums of finite progressions of successive integers raised to
fixed positive powers.  Equally valuable is online Reference [11], which cites a rich literature and
covers besides a vast panorama of diverse mathematical knowledge.

Bernoulli identity (5), which is the principal object of our present concern, springs into view by setting to
zero the imaginary part of the analytic quadrature (6), below, around contour $\,\Omega\,,$
with the corresponding null value requirement on its real part providing an evaluation of the general term from sequence (4),
listed in (14).  No claim whatsoever is made here as to any ultimate novelty in outcome (14), which is available
in symbolic form at any desired index $\,n\,$ through routine demand from {\em{Mathematica}}.  Outcome (14),
expressed here as a finite sum of Riemann zeta functions at odd integer arguments, continues to attract the attention
of contemporary research focused upon polylogarithms [4, 5, 8, 9].  But the formulae thus
made available are subordinated in [4, 5] and elsewhere to the task of evaluating a variety of dissimilar
quantities, and appear to be tangled in thickets of notation.  From this
standpoint, formula (14) (and its identical twin (51) derived in an even more elementary fashion)
may perhaps still provide the modest service of a stand-alone, encapsulated result, easily derived and
easily surveyed.  In particular, the canonical method of derivation evolved in [8] and repeatedly
alluded to in [4, 5] requires rather strenuous differentiations of Gamma function ratios,
and results finally in a recurrence on the individual $\,I_{n}\,$ (or else an equivalent generating function).
To be sure, while the work in [8] is immensely elegant, it is at the same time immensely more intricate
than either of our independent derivations culminating in (14) and (51).

    On the other hand, it does appear to have escaped previous notice that the Bernoulli recurrence (5), which
is ancient and foundational in its own right, should likewise evolve (via (16)) from the same
quadrature around contour $\,\Omega\,$ when one insists that the
corresponding imaginary part also vanish.  And, just as is the case with (14), formula (16), too, is shaped by
its contact with Riemann zeta functions, but evaluated this time at even integer arguments, which latter
circumstance, by virtue of the celebrated Euler connection, opens the portal to entry by the similarly
indexed Bernoulli numbers.  It is of course none of our purpose here to compete with, let alone to supplant
in any way the standard derivations of (5).  Rather, we seek merely to highlight its reappearance in what surely
must be conceded to be an unexpected setting.

   In Section 4 we augment the discussion by admitting the full-fledged Ber-
noulli polynomials $B_{n}(z)$ as
integrand factors during contour integration (Eq. (26)).  And, while this route will no longer lead us directly to
recurrence (5), it will underwrite two key ingredients upon which its demonstration in Sections 2 and 3
pivots, to wit, an {\textit{ab initio}} derivation of the Euler link (17) between Riemann's zeta $\zeta(2m)$
and Bernoulli number $B_{2m}$ at all even indices $2m\geq 2,$ and the odd index nullity $B_{2m+1}=0$
$\,\forall\,m\geq 1$ invoked during passage from Eq. (19) to Eq. (20) below.  Moreover, the toolkit of
Bernoulli polynomial identities will provide, in Eq. (43), a fresh derivation of great elegance, as if
one were still needed, of the log-sine evaluation (14).   

   We round out this note with an appendix wherein contour integration cedes place to the more elementary
setting of a Fourier series on whose basis (14) is recovered yet again (as Eq. (51)) through repeated integration by
parts.  That same Fourier series provides moreover an exceedingly short and simple confirmation of (1),
complementary to the contour integral method, an option to which allusion has already been made in Footnote 3.
Of course, at this point, no further light can, nor need be shed upon (5) \textit{per se}.

\section{Null Quadratures on Contour $\,\Omega$.}

       Guided by the cited example 5 in [1, Section 5.3], we consider for $\,n\,\geq\,0\,$ the sequence of numbers
\begin{equation}
K_{n}\,=\,\int_{\,\Omega}z^{n}\,{\rm{Log}}\left(1-e^{\,2iz}\right)dz\,=\,0\;,
\end{equation}
all of them annulled by virtue of closed contour $\,\Omega\,$ being required to
lie within a domain of analyticity for $\,{\rm{Log}}\left(1-e^{\,2iz}\right)\,$
in the plane of complex $z=x+iy.$  Save for quarter-circle indentations of vanishing radius
$\delta\,$ around $\,z\,=\,0\,$ and $\,z\,=\,\pi\,,$ contour $\,\Omega\,$
bounds a semi-infinite vertical strip, with a left leg having $\,x\,=\,0\,$
fixed and descending from $\,y\,=\,\infty\,$ to $\,y\,=\,\delta\,$ (quadrature
contribution $\,L_{n}\,$), and a right leg at a fixed $\,x\,=\,\pi\,$ ascending from
$\,y\,=\,\delta\,$ to $\,y\,=\,\infty\,$ (quadrature contribution $R_{n}$),
linked at their bottom by a horizontal segment with $\,y\,=\,0\,$ and
$\,\delta\,\leq\,x\,\leq\,\pi-\,\delta\,$ (quadrature contribution $H_{n}$).
In what follows it will be readily apparent that the limit $\,\delta\!\downarrow\!0+$
may be enforced with full impunity, a gesture whose {\textit{fait accompli}} status will
be taken for granted.  Likewise passed over without additional
comment will be the fact that no contribution is to be sought from contour completion
by a retrograde horizontal segment $\,\pi\,\geq\,x\,\geq\,0\,$ at
infinite remove, $\,y\,\rightarrow\,\infty\,.$

     We now find
\begin{equation}
L_{n}\,=\,-\,i^{n+1}\int_{\,0}^{\,\infty}y^{\,n}\log\left(1-e^{-2y}\right)dy\;,
\end{equation}
\begin{equation}
R_{\,n}\,=\,+\,i\int_{\,0}^{\,\infty}\left(\pi+iy\right)^{n}\log\left(1-e^{-2y}\right)dy\;,
\end{equation}
and
\begin{eqnarray}
H_{n} \rule{-1.5mm}{0mm}   &   =  & \rule{-1.5mm}{0mm} \int_{\,0}^{\,\pi}x^{n}\left[\rule{0mm}{4mm}\log(2)-
                         \frac{i\pi}{2}+ix+\log(\,\sin(x)\,)\right]dx      \nonumber     \\
      \rule{-1.5mm}{0mm}   &   =  &  \rule{-1.5mm}{0mm}\frac{\pi^{n+1}}{\,n+1\,}\log(2)-i\frac{\pi^{n+2}}{\,2(n+1)\,}+
i\frac{\pi^{n+2}}{\,n+2\,}+\int_{\,0}^{\,\pi}x^{n}\log(\,\sin(x)\,)dx\,.     
\end{eqnarray}
Series expansion of the logarithm further gives
\begin{equation}
L_{n}\, = \, +\,i^{n+1}\sum_{l\,=\,1}^{\infty}\frac{\,1\,}{\,l\,}\int_{\,0}^{\,\infty}y^{\,n}e^{-2ly}\,dy \,=\,+\,
i^{n+1}\,\frac{\,n!\,}{\,2^{n+1}\,}\sum_{l\,=\,1}^{\infty}\frac{\,1\,}{\,l^{n+2}\,} \,,
\end{equation}
the interchange in summation and integration being legitimated by Beppo Levi's monotone convergence theorem, and similarly
\begin{equation}
R_{\,n}\, = \, -\,i\sum_{k\,=\,0}^{n}\left(\!\!\begin{array}{c}
                                                 n \\
                                                 k
                                         \end{array}\!\!\right)\pi^{n-k}\,i^{k}
\frac{\,k!\,}{\,2^{k+1}\,}\sum_{l\,=\,1}^{\infty}\frac{\,1\,}{\,l^{k+2}\,}\,,
\end{equation}
in both of which there insinuates itself the Riemann zeta function 
\begin{equation}
\zeta\,(s)\,=\,\sum_{m\,=\,1}^{\infty}\frac{\,1\,}{\,m^{s}\,}
\end{equation}
at a variety of its argument values $\,s.$\footnote{This canonical
definition implies a guarantee of series convergence, assured by the requirement that $\,\Re\,s\,>\,1\,.$
A robust arsenal of knowledge exists for continuing $\,\zeta(s)\,$ across the entire plane of
complex variable $\,s\,=\,\sigma+it\,,$ with a simple pole intruding at $\,s\,=\,1\,.$}
So armed, we proceed next to set
\begin{equation}
K_{n}\,=\,L_{n}+H_{n}+R_{n}\,=\,0
\end{equation}
and remark that, regardless of the parity of index $\,n\,,$ $\,L_{n}\,$ {\textit{per se}}
is always absorbed by the contribution from the highest power $\,y^{\,n}\,$ within the
integrand for $\,R_{\,n}\,.$  This circumstance accounts for the imminent appearance of
the floor function affecting the highest value
of summation index $\,k\,$ in Eqs. (14)-(16) and (19) below.

      A requirement that the real part of (13) vanish provides now the following string of valuable
log-sine quadrature formulae
\begin{eqnarray}
\int_{\,0}^{\,\pi}x^{n}\log(\,\sin(x)\,)\,dx     &    =  &    -\frac{\pi^{n+1}}{\,n+1\,}\log(2)         \nonumber   \\
        &      &  \rule{-2.3cm}{0mm}   +\,\frac{n!}{\,2^{n+1}\,}\!\sum_{\,k\,=\,1}^{\lfloor n/2 \rfloor}(-1)^{k}\,  
                                          \frac{(2\pi)^{n-2k+1}}{\,(n-2k+1)!\,}\,\zeta(2k+1)  \;,
\end{eqnarray}
of which the first two, at $\,n\,=\,0\,$ and $\,n\,=\,1\,,$ with the sum on the
right missing, validate (2) and (3), each one of them being in any event widely tabulated.  And again, as was
first stated in Footnote 3, Eq. (14) is consistently reaffirmed by {\em{Mathematica}},
even when harnessed in its symbolic mode.  We note in passing the self-evident fact that,
unlike the corresponding prescriptions found in [\textbf{9,10}],
formula (14) is fully explicit, needing to rely neither upon a generating function nor
a recurrence, even though, naturally, such recurrence arrives at a final rendezvous with identically
the same result.

     A close prelude to identity (5) follows next from the co\"{e}xisting  requirement that
the imaginary part of (13) vanish.  This requirement takes the initial form
\begin{equation}
-\frac{\pi^{n+2}}{\,2(n+1)\,}+\frac{\pi^{n+2}}{\,n+2\,} =
 \sum_{k=0}^{\lfloor \frac{n-1}{2} \rfloor }\left(\!\!\begin{array}{c}
                                                 n \\
                                                 2k
                                         \end{array}\!\!\right)\pi^{n-2k}(-1)^{k}
\frac{(2k)!}{\,2^{2k+1}\,}\,\zeta(2k+2) \rule{0.0cm}{0mm}       
\end{equation}
and is subsequently moulded into the shape
\begin{equation}
\sum_{k=0}^{\lfloor \frac{n-1}{2} \rfloor }\left(\!\!\begin{array}{c}
                                                 n \\
                                                 2k
                                         \end{array}\!\!\right)  \frac{B_{2k+2}}{\,(k+1)(2k+1)\,}\,=\,\frac{n}{\,(n+1)(n+2)\,}    
\end{equation}
on taking note of Euler's
celebrated connection [2, 3, 7, 10, 11, 12]
\begin{equation}
\zeta(2k)\,=\,(-1)^{k+1}(2\pi)^{2k}\frac{B_{2k}}{\,2(2k)!\,}    \;\;\;\forall\, k\geq 1
\end{equation}
allowing us to displace attention from the even-argument values of Riemann's zeta
to the correspondingly indexed Bernoulli numbers $\,B_{2k}\,.$

\section{Recurrence Reduction.}

     Recurrence (16) is not quite yet in the desired form (5), but it is easily steered
toward this goal.  That process begins by noting that
\begin{equation}
\left(\!\!\begin{array}{c}
                    n \\
                   2k
          \end{array}\!\!\right)
\frac{1}{\,(k+1)(2k+1)\,}\,=\,
\left(\!\!\begin{array}{c}
            \rule{1.5mm}{0mm} n+2 \\
                 2k+2
          \end{array}\!\!\right)\frac{2}{\,(n+1)(n+2)\,}\;,
\end{equation}
whereupon (16) becomes
\begin{equation}
\sum_{k\,=\,0}^{\lfloor \frac{n-1}{2} \rfloor}\left(\!\!\begin{array}{c}
                                                   \rule{1.5mm}{0mm} n+2 \\
                                                              2k+2
                                                        \end{array}\!\!\right)B_{2k\,+\,2}\,=\,\frac{\,n\,}{\,2\,}\;.
\end{equation}
Now the advance of $\,2k\,$ in steps of two means that it reaches a maximum value $\,M=n-1\,$ when $\,n\,$ is odd,
and one offset instead by two below $n,$ $M=n-2,$ when $\,n\,$ is even.  At the same time the accepted {\mbox{null value of odd-index}}
Bernoulli numbers starting with $\,B_{3}=0\,$ means that we are free, and self-consistently so, to intercalate all odd indices
missing from the progression $2k+2$ in order to attain an index advance in steps of one
and to entertain a common index maximum of $n+1,$ regardless of the parity of $n.$  Altogether then, (19) admits a restatement as
\begin{equation}
\sum_{k\,=\,2}^{n+1}\left(\!\!\!\begin{array}{c}
                                   n+2 \\
                                    k
                              \end{array}\!\!\!\right)B_{k}\,=\,\frac{\,n\,}{\,2\,}\;,
\end{equation}
or else
\begin{equation}
\sum_{k\,=\,0}^{n+1}\left(\!\!\!\begin{array}{c}
                                   n+2 \\
                                    k
                              \end{array}\!\!\!\right)B_{k}\,=\,\frac{\,n\,}{\,2\,}\,+\,\left\{\left(\!\!\!\begin{array}{c}
                                   n+2 \\
                 \rule{0.5mm}{0mm}  0
                              \end{array}\!\!\!\right)B_{0}\,+\,\left(\!\!\!\begin{array}{c}
                                   n+2 \\
                 \rule{0.5mm}{0mm}  1
                              \end{array}\!\!\!\right)B_{1}\right\}\;.
\end{equation}
But now we find that
\begin{equation}
\left(\!\!\!\begin{array}{c}
                 n+2 \\
\rule{0.5mm}{0mm} 0
          \end{array}\!\!\right)B_{0}\,+\,\left(\!\!\!\begin{array}{c}
                 n+2 \\
\rule{0.5mm}{0mm} 1
          \end{array}\!\!\right)B_{1}\,=\,1\,-\,\frac{\,n+2\,}{\,2\,}\,=\,-\frac{\,n\,}{\,2\,}\,,
\end{equation}
with the effect of reducing (21) to just
\begin{equation}
\sum_{k\,=\,0}^{n+1}\left(\!\!\!\begin{array}{c}
                                   n+2 \\
                                    k
                              \end{array}\!\!\right)B_{k}\,=\,0\;,
\end{equation}
which is nothing other than (5).

\section{Bernoulli Polynomials as Integrand Factors on Contour $\,{\tilde{\Omega}}$.}

    Added perspective accrues when power multiplier $z^{n}$ in the integrand from (6) is replaced by the Bernoulli polynomial
$B_{n}(z).$  By following this route we will, on the one hand, forfeit a direct access to recurrence (5), but, by way of
compensation, we will recover Euler's connection (17) linking even index Bernoulli number $B_{2m}$ values to their
Riemann $\zeta(2m)$ counterparts, and will similarly confirm the vanishing $B_{2m+1}=0$ of all odd index Bernoulli numbers
beginning with index 3.  More even than that, a continuing interplay between contour integration around a half strip as presently
considered, and one among the tangle of Bernoulli polynomial identities, will disclose once more the log-sine formula (14).

      All required Bernoulli polynomial relations will be drawn without further comment from [7], which is readily accessed electronically.
And so, bowing to
the notational convention adopted there, we demote the Bernoulli number symbol from capital $B_{n}$ to lower case $b_{n},$
upper case $B_{n}(z)$ being reserved for the polynomials so named.\footnote{This shift evidently seeks to minimize visual confusion.}
The polynomials themselves are uniquely defined by
\begin{equation}
\left\{
\begin{array}{ll}
B_{0}(z)=1                                                   \\
%                                                            \\
 B_{n}(z+1)-B_{n}(z)=nz^{n-1},       & \forall\, n \geq 1    \\
%                                                            \\
 {\Large{\int_{0}^{1}}}B_{n}(t)dt=0, &  \forall\, n \geq 1
\end{array}
\right.
\end{equation}
[7, Section 1, Corollary 1.3] and the Bernoulli numbers $b_{n}$ then follow from
\begin{equation}
b_{n}=B_{n}(0)\;\;\forall\, n\geq 0
\end{equation}
[7, Section 2, Definition 2.2].  Their recurrence (5) is found as entry {\textit{iv}}
to [7, Section 2, Proposition 2.3], with no mention whatsoever at that point of a Riemann
$\zeta$ connection.  Also disclosed there under entry {\textit{i}} is the vanishing of all odd index numbers $b_{2n+1}=0$
$\,\forall\, n \geq 1.$  The auxiliary fact that $\Im(B_{n}(z))=0$ and hence {\em{a fortiori}} that $\Im(b_{n})=0$ without
exception whenever $\Im(z) = 0,$
follows from the bijective mapping established between real polynomial bases in [7, Section 1], Lemma 1.1 and
Bernoulli polynomial Definition 1.2, which latter underwrites the constructive recipe in (24).

     Our imminent intent to capitalize on the second entry in (24) suggests that width $\pi$ of
contour $\Omega$ be compressed to just $1$ under variable scaling $x=\pi u,$ $y=\pi v,$ acknowledged by the
notational shift $\Omega\rightarrow{\tilde{\Omega}}$ in the plane of complex $\xi=u+iv.$  And so
we replace quantities $K_{n}$ from (6) by similarly null analogs
\begin{equation}
{\tilde{K}}_{n}\,=\,\int_{\,{\tilde{\Omega}}}B_{n}(\xi)\,{\rm{Log}}\left(1-e^{2i\pi\xi}\right)d\xi\,=\,0\,.
\end{equation} 
The contour integration apparatus of Eqs. (7)-(13), marshalled out on behalf of $K_{n},$ carries over essentially verbatim, ceding place now to
\begin{equation}
{\tilde{L}}_{n}\,=\,-\,i\int_{\,0}^{\,\infty}B_{n}(iv)\log\left(1-e^{-2\pi v}\right)dv\;\,,
\end{equation}
\begin{equation}
{\tilde{R}}_{\,n}\,=\,+\,i\int_{\,0}^{\,\infty}B_{n}(1+iv)\log\left(1-e^{-2\pi v}\right)dv\;\,,
\end{equation}
and
\begin{eqnarray}
{\tilde{H}}_{n}    &   =  &  \int_{\,0}^{\,1}B_{n}(u)\left[\rule{0mm}{4mm}\log(2)-
                         \frac{i\pi}{2}+i\pi u+\log(\,\sin(\pi u)\,)\right]du     \nonumber  \\
                   &   =  &  \int_{\,0}^{\,1}B_{n}(u)\left[\rule{0mm}{4mm}\log(2)+\log(\,\sin(\pi u)\,)+i\pi B_{1}(u)\right]du 
\end{eqnarray}
wherein we have identified $(u-1/2)$ with $B_{1}(u)$ [7, Section 1, Corollary 1.3].  And, of course, we still have, $\,\forall\, n\geq 0,$
\begin{equation}
{\tilde{K}}_{n}\,=\,{\tilde{L}}_{n}+{\tilde{H}}_{n}+{\tilde{R}}_{n}\,=\,0
\end{equation}
which, in particular, forces one to consider the sum\footnote{Result (31) holds only for $n\geq 1.$  When $n=0$ the sum
${\tilde{L}}_{0}+{\tilde{R}}_{\,0}$ is trivially null by virtue of the first line in (24).  The parent quantities ${\tilde{K}}_{n},$
by contrast, are null $\,\forall\, n\geq 0$ without restriction.}
\begin{eqnarray}
{\tilde{L}}_{n}+{\tilde{R}}_{\,n}& = &  i\int_{0}^{\infty}\left\{\rule{0mm}{4mm}B_{n}(1+iv)-B_{n}(iv)\right\}\log\left(1-e^{-2\pi v}\right)dv \nonumber \\
                                 & = &  ni^{n}\int_{0}^{\infty}v^{n-1}\log\left(1-e^{-2\pi v}\right)dv \nonumber \\
                                 & = &  \frac{ni^{n}}{(2\pi)^{n}}\int_{0}^{\infty}v^{n-1}\log\left(1-e^{-v}\right)dv  \nonumber \\
                                 & = & -\frac{ni^{n}}{(2\pi)^{n}}\sum_{l=1}^{\infty}\frac{1}{\,l\,}\int_{0}^{\infty}v^{n-1}e^{-lv}dv   \\
                                 & = & -\frac{n!i^{n}}{(2\pi)^{n}}\sum_{l=1}^{\infty}\frac{1}{\,l^{n+1}\,}  \nonumber \\
                                 & = & -\frac{n!i^{n}}{(2\pi)^{n}}\zeta(n+1)  \nonumber
\end{eqnarray}
after invoking the second line from (24) and mimicking the series expansion in (10).

      Turning attention once more to (29), we must again discriminate between the null index $n=0$ case and all others with $n\geq 1.$  Thus, when $n=0,$
reference to the first and third lines in (24) and Footnote 6 gives
\begin{equation}
{\tilde{K}}_{0}={\tilde{H}}_{0}=\log(2)+\int_{0}^{1}\log(\,\sin(\pi u)\,)\,du=0
\end{equation}
or else
\begin{equation}
\int_{0}^{\pi}\log(\,\sin(x)\,)\,dx=-\pi\log(2)\,,
\end{equation}
which is Eq. (2).

      By contrast, for $n\geq 1,$ recourse to [7, Section 3, Corollary 3.3] provides, when integer indices $p$ and $q$ are both in excess of zero\footnote{Both
sides of (34) are symmetric in indices $p$ and $q,$ on its left by inspection while on its right because of the fact, soon to be confirmed yet again (Eqs. (38) and (40) below),
that $b_{p+q}$ vanishes unless $p+q$ is even.}  
\begin{equation}
\int_{0}^{1}B_{p}(u)B_{q}(u)\,du=\frac{(-1)^{q-1}}{\left(\!\!\!\begin{array}{c}
                                                              p+q \\
                                                               q
                                                         \end{array}\!\!\!\right)}\,b_{p+q}
\end{equation}
\parindent=0.0in
and thus causes (29) to read
\begin{equation}
{\tilde{H}}_{n}=\int_{0}^{1}B_{n}(u)\log(\,\sin(\pi u)\,)\,du+\frac{i\pi}{\,n+1\,}\,b_{n+1}\,.
\end{equation}
And so, on putting (30), (31), and (35) together we arrive at
\begin{equation}
\int_{0}^{1}B_{n}(u)\log(\,\sin(\pi u)\,)\,du+\frac{i\pi}{\,n+1\,}\,b_{n+1}=\frac{n!i^{n}}{\,(2\pi)^{n}\,}\,\zeta(n+1)
\end{equation}
whenever $n\geq 1.$  On sorting (36) out according to its real and imaginary components and the parity of $n$, we thus find,
when $n=2m,$
\begin{eqnarray}
\int_{0}^{1}B_{2m}(u)\log(\,\sin(\pi u)\,)\,du &  =  &  \frac{\,(2m)!(-1)^{m}\,}{\;\;(2\pi)^{2m}}\,\zeta(2m+1)  \\
b_{2m+1} & = & 0   
\end{eqnarray}
and if instead $n=2m-1,$
\begin{eqnarray}
\int_{0}^{1}B_{2m-1}(u)\log(\,\sin(\pi u)\,)\,du &  = & 0    \\
b_{2m}   & = & \frac{\,2(2m)!(-1)^{m+1}}{\;\;(2\pi)^{2m}}\,\zeta(2m)  
\end{eqnarray}
$\forall\, m\geq 1.$  Odd index vanishing of $b_{2m+1}$, $m\geq 1,$ is thus vindicated, as is also the full content of Euler's connection (17)
at all even indices $b_{2m}$ beginning with $m=1.$\footnote{The missing start-up values $b_{0}=1$ and $b_{1}=-1/2$ follow respectively
from $b_{0}=B_{0}(0)=1$ and $b_{1}=B_{1}(0)=(u-1/2)|_{u=0}=-1/2.$}    
\parindent=0.20in

     There exists one further, polynomial basis inversion identity, namely
\begin{equation}
u^{n}=\frac{1}{n+1}\sum_{m=0}^{n}\left(\!\!\!\begin{array}{c}
                                                  n+1 \\
                                                   m
                                             \end{array}\!\!\!\right)B_{m}(u)   \;\;\; \forall\, n\geq 0
\end{equation}
[7, Section 2, Application 1], which enables us to consolidate the fragmented information
scattered among the quadratures (37) into a single, compact form.  Thus, with (32) adjoined,
\begin{eqnarray}
\int_{0}^{1}u^{n}\log(\,\sin(\pi u)\,)\,du &  =  &  -\frac{1}{\,n+1\,}\log(2)             \nonumber   \\
          &      & \rule{-2.3cm}{0mm} +\,n\,!\sum_{m=1}^{\lfloor n/2 \rfloor}\frac{(-1)^{m}}{\,(n-2m+1)!(2\pi)^{2m}\,}\,\zeta(2m+1)  
\end{eqnarray}
or else
\begin{eqnarray}
\int_{0}^{\pi}x^{n}\log(\,\sin(x)\,)\,dx  &  =  & -\frac{\pi^{n+1}}{\,n+1\,}\log(2)         \nonumber   \\
              &        &  \rule{-2.3cm}{0mm} + \frac{n!}{2^{n+1}}\sum_{m=1}^{\lfloor n/2 \rfloor}(-1)^{m}\frac{(2\pi)^{n-2m+1}}{\,(n-2m+1)!\,}\,\zeta(2m+1)  
\end{eqnarray}
which is (14) once more.

      The null quadratures in (39) are a consequence of the antisymmetry
\begin{equation}
B_{2m-1}(1/2-\Delta u)=-B_{2m-1}(1/2+\Delta u)
\end{equation}
around $u=1/2$ which follows readily from the fact that\footnote{Eq. (44) is clearly accompanied by a symmetric counterpart
around $u=1/2$ for the even-indexed polynomials.  These symmetry/antisymmetry attributes are corroborated under a different guise
by the Fourier series evolved in [7, Section 3, Proposition 3.1, parts {\textit{i}} \& {\textit{ii}}].
Consistent too with (44) are the null
mid-point evaluations $B_{2m-1}(1/2)=(4^{1-m}-1)b_{2m-1}=0$ $\,\forall\, m \geq 1$ [7, Section 2, Proposition 2.3{\textit{ii}}],
with the special evaluation $B_{1}(1/2)=0$ being already self-evident from $B_{1}(u)=u-1/2.$  For the even-indexed polynomials, the
overarching, global constraint $\int_{0}^{1}B_{2n}(u)du=0$ $\,\forall\, n\geq 1$ is maintained despite their exhibiting non-zero
mid-point values $B_{2n}(1/2)=(2^{1-2n}-1)b_{2n}\neq 0\;\;\forall\, n\geq 0.$}
\begin{equation}
B_{n}(1-u)=(-1)^{n}B_{n}(u)\;\;\forall\,n\geq 0
\end{equation}
[7, Section 2, Proposition 2.1{\textit{ii}} (augmented so as to include $n=0$)]. 
When viewed at $m=1$ and with reference to (32), Eq. (39) leads us thus to confront
\begin{eqnarray}
\int_{0}^{1}u\log(\,\sin(\pi u)\,)\,du-\frac{1}{2}\int_{0}^{1}\log(\,\sin(\pi u)\,)\,du    &     &     \nonumber   \\
        &        & \rule{-4.7cm}{0mm}    =  \, \int_{0}^{1}u\log(\,\sin(\pi u)\,)\,du+\frac{1}{2}\log(2) \, = \, 0
\end{eqnarray}
which amounts to
\begin{equation}
\int_{0}^{\pi}x\log(\,\sin(x)\,)\,dx=-\frac{\pi^{2}}{2}\log(2)
\end{equation}
in agreement with both (3) and (43), the trailing sum of the latter being then vacuous.  Precursor (33) is,
{\textit{ipso facto}}, likewise subsumed under (43) by the very manner of the latter's derivation.

      And finally, even though it has now been sidestepped, recurrence (5) is duly assembled in
[7, Section 2, Proposition 2.3{\textit{iv}}],
free from any reference to Riemann's $\zeta.$  The basis of deduction at that point is still another Bernoulli polynomial identity
\[
B_{n}(u)=\sum_{m=0}^{n}\left(\!\!\!\begin{array}{c}
                                                   n \\
                                                   m
                                             \end{array}\!\!\!\right)b_{n-m}u^{m}   \;\;\; \forall\, n\geq 0  \;,
\]
reinforced by the observation, also derived there, that $B_{n}(1)=b_{n}$ $\,\forall\, n \geq 2.$
Verily, verily, on this arena all roads seem to converge upon the proverbial Rome.
\vspace{-0.5mm}

\section{Appendix:  A Fourier Series Grace Note.}

    A somewhat more pedestrian derivation of (14) rests upon consideration of the power series
\begin{equation}
{\rm{Log}}\,(\,1-z\,)\,=\,-\,\sum_{l\,=\,1}^{\infty}\,\frac{\,z^{\,l}\,}{\,l\,}
\end{equation}
along the unit circle $\,z\,=\,e^{\,i\,\vartheta}.$
Separation into real and imaginary parts resolves itself into a pair of Fourier series
\begin{equation}
\log\left\{\rule{0mm}{4mm}2\left|\rule{0mm}{3mm}\sin(\vartheta/2)\right|\right\}\,=\,
-\sum_{l\,=\,1}^{\infty}\,\frac{\,\cos(l\vartheta)\,}{\,l\,}
\end{equation}
and
\begin{equation}
\frac{\,1\,}{\,2\,}\left\{\rule{0mm}{4mm}\vartheta\,({\rm{mod}}\,2\pi)-\pi\right\}
=\,-\sum_{l\,=\,1}^{\infty}\,\frac{\,\sin(l\vartheta)\,}{\,l\,}\;\,,
\end{equation}
of which the second is of no interest {\textit{vis-\`{a}-vis}} our immediate objective.
The logarithmic divergence on both left and right in (49) whenever $\,\vartheta\,\equiv\,0\,({\rm{mod}}\;2\pi)$
remains integrable and is thus taken henceforth in easy stride.

      Repeated integration by parts {\textit{vis-\`{a}-vis}} series (49),
when first multiplied by the argument power $\,\vartheta^{\,n},$
advances by $\,\cos\rightarrow\sin\rightarrow\cos\,$ couplets, with end-point contributions
arising only on the second beat, and the argument powers falling
in steps of two.\footnote{In particular, this quadrature cadence provides a
motivation, alternative to that previously given,
as to why it is that the floor function affects the upper index cutoff
$\,\lfloor n/2 \rfloor\,$ in both (14) and (51), allowing for unit growth in that
cutoff only when $\,n\,$ {\textit{per se}} advances by two.}  One assembles in this manner the general
formula
\begin{eqnarray}
\int_{\,0}^{\,\pi}\,\vartheta^{\,n}\log(\,\sin(\vartheta)\,)\,d\,\vartheta & = & -\frac{\,\pi^{n+1}\,}{\,n+1\,}\log(2)   \nonumber    \\
          &      & \rule{-2.5cm}{0mm}+\frac{n!}{2^{n+1}}\sum_{k=1}^{\lfloor n/2 \rfloor}(-1)^{k}  
                       \frac{\,(2\pi)^{n-2k+1}\,}{\,(n-2k+1)!\,}\,\zeta(2k+1) 
\end{eqnarray}
holding good unrestrictedly for $\,n\,$ even or odd, and agreeing in every respect with (14).  The only
wrinkle to notice, perhaps, is that the sequence of integrations by parts which underlies (51) terminates,
at each summation index $\,l\,$ in (49), with a term proportional to either
\begin{equation}
\int_{\,0}^{\,\pi}\cos(\,2l\vartheta\,)\,d\,\vartheta\,=\,0
\end{equation}
in the event that $\,n\,$ is even, or
\begin{equation}
\int_{\,0}^{\,\pi}\vartheta \cos(\,2l\vartheta\,)\,d\,\vartheta\,=\,0
\end{equation}
otherwise.  Equation (52) is of course obvious whereas (53), while equally true and
welcome as such, is, at first blush, mildly surprising.  All in all the derivation which underlies (14) is
far smoother and less apt to inflict bookkeeping stress, even if it is (51) which seems to rest
on a more elementary underpinning.

    It would be truly disappointing were we not able to utilize (49) so as to give an
essentially one-line, zinger-style proof of (1).  This anticipation is readily met simply
by squaring both sides of (49), with
summation indices $\,l\,$ and $\,l\,'\,$ figuring
now on its right, and noting that when, as here, both $\,l\,\geq\,1\,$ and
$\,l\,'\,\geq\,1\,,$
\begin{equation}
\int_{\,0}^{\,\pi}\cos(\,2l\vartheta\,)\cos(\,2l'\vartheta\,)\,d\vartheta=\frac{\pi}{2}\delta^{\,l}_{\,l'}\;\,,
\end{equation}
with $\,\delta^{\,l}_{\,l'}\,$ being the Kronecker delta, unity when its indices match, and zero otherwise.
In a gesture which embodies the essence of Parseval's theorem, it follows immediately that
\begin{equation}
\frac{\,1\,}{\,2\,}
\int_{\,0}^{\,\pi}\left\{\rule{0mm}{4mm}\log(\,2\,\sin(\vartheta)\,)\right\}^{\,2}\!d\vartheta\,=\,
\frac{\,\pi\,}{\,4\,}\sum_{l=1}^{\infty}\frac{\,1\,}{\,l^{\,2}\,}\,=\,\frac{\;\;\pi^{\,3}\,}{\,24\,}\;\,,
\end{equation}
and we are done.
\parindent=0.0in
\vspace{6mm}

{\textbf{Acknowledgement.}}  Thanks are due to an anonymous referee for suggesting that our viewpoint would benefit
from having the Bernoulli polynomials brought into play.

\end{document}